\documentclass[12pt]{article}
\usepackage{latexsym}
\usepackage{amsfonts}
\usepackage{amssymb}
\newtheorem{theorem}{Theorem.}[section]

\begin{document}

\title{Some non-finitely presented Lie Algebras}
 
\author{Joseph Abarbanel and Shmuel Rosset\\ 
Tel-Aviv University\\
Ramat-Aviv 69978\\
Israel\\
email: rosset@math.tau.ac.il,\\
yossia@math.tau.ac.il}
\date{}
 
\maketitle

\abstract

  Let $L$ be a free Lie algebra over a field $k$, $I$ a non-trivial 
 proper ideal of $L$, $n>1$ an integer. The multiplicator
 $H_2(L/I^n,k)$ of $L/I^n$ is not finitely generated, and so in
 particular, $L/I^n$ is not finitely presented, even when $L/I$
 is finite dimensional.

\section {Introduction}

If $R$ is a free associative algebra, over a field, and $I$ is a two 
sided ideal of $R$, then Lewin proved \cite{Le} that $I^2$ is not 
finitely generated (as a 2-sided ideal!) when the algebra $R/I$ 
is infinite dimensional. In other words, $R/I^2$ is not finitely 
presented in this case. On the other hand, it is easy
to see that when $R$ is finitely generated and $R/I$ is finite dimensional,
so is $R/I^2$, and hence $I^2$ is finitely generated.

 Similar behavior is seen in groups. If $F$ is a finitely generated
free group, and $R$ is a normal subgroup then $R'$ is {\em normally}
finitely generated if, and only if, $F/R$ is finite. In fact Baumslag,
Strebel and Thomson proved \cite{Ba} a stronger fact. Denoting the
m-th member of the lower central series by $\gamma_m$, they proved that for
$m>1$ the Schur multiplier of $F/\gamma_mR$, $H_2(F/\gamma_mR,\mathbb Z)$, 
is notfinitely generated (as an abelian group) if $F/R$ is not finite.
\\
We note that for the three statements

(a) $R$ is normally finitely generated,

(b) $R/R'$ is finitely generated as a module over $G=F/R$,

(c) $H_2(G,\mathbb Z)$ is finitely generated as an abelian group\\
we have $(a)\Rightarrow (b)\Rightarrow (c)$.

 In this paper we prove a result of similar nature for Lie algebras. 

\begin{theorem}
\label{th:main}
Let $L$ be a free Lie algebra with basis $X$, over a field $k$, and $I$ be {\em any} 
non-zero proper ideal of $L$, then $I'=[I,I]$ is not finitely generated 
as an ideal. In fact, the ``Schur multiplier'' of $L/I^n$, 
$H_2(L/I^n,k)$, is not finitely generated if $n>1$, and hence $L/I^n$ is
not finitely presented.
\end{theorem}
Here $I^n$ denotes $I$, if $n=1$, and $[I^{n-1},I]$
if $n>1$. Our proof closely follows the lines of \cite{Ba}.

  In \S \ref{sec:prems} we define some notations and the Magnus
embedding. In \S \ref{sec:image} we build a mapping from the Schur
multiplier into a tensor product of $n-1$ copies of $U(L/I)$. This is
similar to the mapping defined in \cite{Ba}. In \S \ref{sec:hopf}
we build a specific isomorphism of Hopf modules, keeping in mind that the
enveloping algebra of a Lie algebra is a Hopf algebra. In \S 
\ref{sec:comp} we employ the mapping and show that the image of the
``Schur multiplicator'' is not finite dimensional, thus proving the
theorem.

  We wish to thank Alon Wasserman for his help in \S \ref{sec:hopf}.

\section {Preliminaries and Notations}
\label{sec:prems}

Let ${\cal G}$ be a Lie algebra. We will denote the Lie multiplication
of two elements $a,b\in {\cal G}$ by $[a,b]$. As we will also be considering
the enveloping algebra of ${\cal G}$, the multiplication in $U({\cal G})$ will 
be denoted simply as $ab$, while the action of an element $l\in U({\cal G})$ on
an element $a\in {\cal G}$ will be denoted by $a\cdot l$. Note that the action is 
the adjoint action, so that if $l\in L$ then
$a\cdot l = [a,l]$.

 Let ${\cal G}$ be a Lie algebra over a field $k$, $U({\cal G})$ its enveloping 
algebra, $\delta U({\cal G})$ the augmentation ideal
of $U$. Suppose $0\to I\to L\to {\cal G}\to 0$ is a free presentation of
${\cal G}$, where $L$ is the free Lie algebra with basis $X$. The
enveloping algebra, $U(L)$, is therefore a free associative algebra, 
with basis $X$, and $\delta U(L)$ is a free $U(L)$ module, 
with a basis in one-to-one correspondence with $X$. Note that over a 
field,if ${\cal G} \not =0$, $U({\cal G})$ is infinite dimensional, and 
is without zero divisors.

 In addition, if ${\cal G}$ is a Lie algebra over a field and $U({\cal G})$ 
is its enveloping algebra, let $U_n({\cal G})$ be the subspace of $U({\cal G})$ 
spanned by all the products of at most $n$ factors from ${\cal G}$. This 
gives a well known ascending filtration of $U({\cal G})$, and we can define 
the {\em degree} of an element $l$ to be the {\em least} integer $n$ 
such that $l\in U_n({\cal G})$. This function has the properties:

1) ${\rm deg}(a+b)\leq max\{{\rm deg}(a),{\rm deg}(b)\}$,

2) if ${\rm deg}(a)<{\rm deg}(b)$ then ${\rm deg}(a+b)={\rm deg}(b)$,

3) ${\rm deg}(ab)={\rm deg}(a)+{\rm deg}(b)$.\\
In particular, if $x\in {\cal G}$ is non-zero then the degree of $x$ is $1$, so if
$x_1,x_2,\ldots,x_n\in {\cal G}$ are all non-zero then 
${\rm deg}(x_1x_2\cdots x_n)=n$.

 Via the adjoint action, $I/I'$ carries the structure of a $U(L)$ module, and
$I$ acts trivially. All modules will be right modules. Therefore 
$I/I'$ is a $U(L/I)$ module in a natural way.
There is a well known embedding of $U(L/I)$ modules, the Magnus embedding,
described below, of $I/I'$ into $\delta U(L)\otimes_{U(L)} U(L/I)$. 
This embedding will be denoted by $\phi:I/I'\to \delta U(L)\otimes_{U(L)} U(L/I)$. 
The action of $L$ on $\delta U(L)\otimes_{U(L)} U(L/I)$ is by right 
multiplication in the right hand term.

 The embedding can be defined in the following way. First define
$\phi:I\to \delta U(L)\otimes_{U(L)} U(L/I)$ by $\phi(x)=x\otimes 1$.
By using the Poincare-Birkhoff-Witt theorem, and the structure it 
gives to $U(L)$, it can be seen that this is a mapping of $U(L)$ modules, i.e.
$\phi(a\cdot l)=\phi(a)l$. First we check the
statement for elements of $L$. If $l\in L$ then $a\cdot l=[a,l]$ and
$\phi([a,l])=[a,l]\otimes 1=(al-la)\otimes 1= a\otimes l - l\otimes a$.
However, $a=0$ in $U(L/I)$ so $\phi([a,l])=a\otimes l=(a\otimes 1)l=\phi(a)l$.
Consider now the subalgebra $A=\{ u\in U(L) | \phi(x\cdot u)=\phi(x)u \quad
 \forall x\in I \}$. Since $L\subset A$ then $A=U(L)$, thus $\phi$ is 
a $U(L)$ module homomorphism.

It is left to show that ${\rm ker} \phi=I'$. If $x\in I'$ then $x$ can be written as $x=\sum [a_i,b_i], a_i,b_i\in I$, so
that $\phi(x)=x\otimes 1= \sum [a_i,b_i] \otimes 1 =
  \sum (a_i b_i - b_i a_i) \otimes 1 = \sum a_i\otimes b_i - b_i \otimes a_i$.
Since $a_i,b_i \in I$ then their images in $U(L/I)$ are $0$ so that 
$\phi(x)=0$.

 Therefore $I' \subset {\rm ker} \phi$. On the other hand suppose 
$x\in {\rm ker} \phi$. Since $\delta U(L)$ is a free $U(L)$ module with 
basis $\{x_i\}$ where $x_i$ is a basis of $L$ as a free Lie algebra, 
we have $x\otimes 1=\sum x_i\otimes f_i$, where, since $\phi(x)=0$, 
$f_i =0$ in $U(L/I)$. Let us denote by $\tilde I$ the kernel of the mapping $U(L)\to U(L/I)$, so that  $f_i\in \tilde I$. But $\tilde I=U(L)I=IU(L)$ and thus
by the Poincare-Birkhoff-Witt theorem this kernel is a 
free left and right $U(L)$ module with a basis that is a basis of 
$I$ as a subalgebra of $L$. Therefore $f_i= \sum w_{i,j} a_j$ where $a_j$ 
are a basis of $I$. It follows that $x=\sum x_i w_{i,j} a_j$.
Consider now the image of $x$, $\bar x$, in $I/I'$. Since $I/I'$ is the commutative Lie
algebra with a basis that is a basis of $I$ as a subalgebra of $L$, then
$\bar x=\sum \lambda_j a_j$, where $\lambda_j\in k$. In other words 
$x=\sum \lambda_j a_j + w,w\in I'$. But
since $I'\subset {\rm ker} \phi$ then we can assume $x=\sum \lambda_j a_j$.
On the other hand $\phi(x)=0$ so $x=\sum x_i w_{i,j} a_j$.
Since $\tilde I$ is a free $U(L)$ module with basis $a_i$
we have $\lambda_j=\sum x_i w_{i,j}$, but $x_i\in \delta U(L)$, so 
$\lambda_j=0$. Hence $x\in I'$, therefore ${\rm ker} \phi = I'$. 

  Another proof of the fact that ${\rm ker}\phi=I'$ can be found in 
\cite{Be} \S 8, as the Magnus embedding is a special case of the 
derivations defined there.

 Throughout the remainder of this paper $I$ will
be a proper non-zero ideal of $L$, and $n>1$ will be an integer.

\section {An image of $H_2(L/I^n,k)$}
\label{sec:image}

 Consider $H_2(L/I^n,k)$. It is known (e.g. \cite{We} p.233) that 
the analogue of the Hopf formula for groups holds for Lie algebras. Therefore
$$ H_2(L/I^n,k) = I^n/[I^n,L] = (I^n/I^{n+1})\otimes_{U(L)} k$$

 We know from the \u Sir\u sov-Witt theorem (see e.g. \cite{Re} p.44) that $I$ 
is a free Lie algebra. Hence $I^n/I^{n+1}$ is, in a natural way, 
identifiable with the $n$-th homogeneous component of the free Lie 
algebra with basis that is a basis of $I/I'$ as a vector space. Since 
the free Lie algebra of a free module can be embedded in the tensor 
algebra over this module, the $n$-th homogeneous component can be 
embedded into the $n$-fold tensor product, i.e. $I^n/I^{n+1}$ can be 
embedded in $\otimes^n I/I'$, where the tensor is over $k$. Any unadorned
tensor product below is to be taken to be over $k$. 
We need this embedding to be a $U(L/I)$ module homomorphism, and it is easy
to see that this is indeed the case when $U(L/I)$ acts on $I^n/I^{n+1}$ via
the adjoint action, and on $\otimes^n I/I'$ diagonally.
The module $\otimes^n I/I'$ can again can be embedded,
through the Magnus embedding, into
$$\otimes^n (\delta U(L)\otimes_{U(L)} U(L/I))$$
Tensoring this with $k$ over $L$ we
get a mapping 
$$ H_2(L/I^n,k) \approx \otimes^nI/I'\otimes_{U(L)} k \to
\otimes^n(\delta U(L)\otimes_{U(L)} U(L/I)) \otimes_{U(L)} k$$

  Since $\delta U(L)$ is a free $U(L)$ module, with a basis $X$ that is a basis
of $L$ as a Lie algebra, we can define for each $x\in X$ a projection, denoted
$p_x : \delta U(L)\otimes_{U(L)} U(L/I) \to U(L/I)$. We therefore have for each
$n$-tuple $(x_1,x_2,\ldots ,x_n)\in X^n$ a mapping
$\phi_{x_1,\ldots ,x_n}:=(p_{x_1}\otimes \cdots \otimes
p_{x_n}\otimes 1)\circ \phi$
$$\phi_{x_1,x_2,\ldots ,x_n}:H_2(L/I^n,k)\to \otimes^nU(L/I) \otimes_{U(L)} k$$

  Since $I/I' \to U(L/I)\otimes \delta U(L)$ is an embedding, there exist
elements $\alpha \in I/I'$ and $x\in X$ such that under the Magnus embedding
and the projection by $x$ the image $a=\phi_x (\alpha)$ is non-zero. These
elements will be put to use below.

\section {Isomorphism of Hopf modules}
\label{sec:hopf}

 As seen in the last section the image of the multiplicator lies
in $(U(L/I)\otimes U(L/I) \cdots \otimes U(L/I))\otimes _{U(L)} k$. On the other
hand it is well known that the enveloping algebra is a Hopf algebra, and
the action with which this module is endowed is consistent with
the standard Hopf structure on $U(L/I)$, which is the diagonal action. 
We shall use the following notation for the structure of Hopf algebras and modules. Let $H$ be a Hopf algebra and $M$ a Hopf module over $H$. 
The diagonal mapping of $H$ will be denoted by $\Delta$, and 
the $n$-fold application of $\Delta$ by $\Delta_n$ (by the
co-associativity of $H$ the components on which we apply $\Delta$ each time
do not matter). The co-unit of $H$ will be denoted by $\epsilon$
(also sometimes known as the augmentation). 
The antipode map of $H$ will be denoted by $S$. The usual
action of $H$ on $M$ will be denoted by multiplication on the right,
and the co-action of $M$ will be denoted by $\rho$. If $h\in H$ then 
$\Delta(h)$ will be written as $\Delta(h)=\sum_{i=1}^l h_1^i\otimes h_2^i$,
and $\Delta(h_1^i)=\sum_{j=1}^{l(i)} h_{1,1}^{i,j}\otimes h_{1,2}^{i,j}$.
If $m\in M$ then $\rho(m)=\sum m_0^i\otimes m_1^i$.

 It is known ( see e.g. \cite{Mo} p.15 ) that for any Hopf algebra $H$ and Hopf
module $M$, $M \approx M'\otimes H$,
where $M'=\{m\in M | \rho(m)=m\otimes 1\}$ with the isomorphism
$m\mapsto \sum m_0^i \cdot S(m_{1,1}^{i,j}) \otimes m_{1,2}^{i,j}$,
where this is actually a double sum on both $i$ and $j$.
It should also be noted that  $M'\otimes H$ is a trivial Hopf
module, i.e. one for which $(m\otimes h)l = m\otimes hl$.

If we now also tensor with $k$ over $H$
we will get 
$$ M \otimes_H k \approx (M' \otimes H)\otimes_H k.$$
However, since $M'\otimes H$ is a trivial (in the sense defined above)
Hopf module we get
$$M\otimes_H k\approx (M'\otimes H)\otimes_H k \approx M'\otimes (H\otimes_H k)\approx M'.$$
The isomorphism is
$$m\otimes 1\mapsto \sum m_0^i \cdot S(m_{1,1}^{i,j})\otimes m_{1,2}^{i,j}\otimes 1\mapsto \sum m_0^i \cdot S(m_{1,1}^{i,j}) \epsilon(m_{1,2}^{i,j}) =$$
$$ \sum m_0^i \cdot S(m_1^i).$$

If we take $M=W\otimes H$ with $W$ any Hopf module, $H$ acting with
the diagonal action and
$$\rho(w\otimes h) = w\otimes \Delta(h)$$
then $M'=W \otimes k\approx W$. 
In this case, if $m=w \otimes h$ then
$\rho(w\otimes h)=w\otimes \Delta(h)$ so
$m_0^i=w \otimes h_1^i$
and $m_1^i=h_2^i$. Therefore
the explicit form of the isomorphism is
$$w\otimes h\otimes 1 \mapsto \sum (w\otimes h_1^i)\Delta(S(h_2^i)).$$

However, we know that the image is in $M'$, so we can apply $1\otimes \epsilon$
to the image and not change it. Also if $h\in H$ then from the 
definition of a Hopf algebra $(1\otimes \epsilon)(\Delta(h))=h\otimes 1$

 Therefore the image is
$$(1\otimes \epsilon)[\sum (w\otimes h_1^i)\Delta(S(h_2^i))]=$$
$$\sum (w\otimes \epsilon(h_1^i))[(1\otimes \epsilon)(\Delta(S(h_2^i)))]=
\sum (w\otimes 1) (\epsilon(h_1^i) S(h_2^i) \otimes 1)= $$
$$(w\otimes 1)(S(h)\otimes 1)$$
so the image in $W$ is
$$ w\otimes h\otimes 1\mapsto wS(h).$$

 In our case we are interested in the module $\otimes^n H$, so we
can take $W=\otimes^{n-1}H$ and the isomorphism will be
$$h_1\otimes h_2 \otimes \cdots \otimes h_n \otimes 1\mapsto
(h_1\otimes h_2 \otimes \cdots \otimes h_{n-1})\Delta_{n-1}(S(h_n)).$$

\section {Computations}
\label{sec:comp}

 We can now prove theorem \ref{th:main}, i.e. show that $H_2(L/I^n,k)$ 
is not finitely generated by exhibiting an infinite number of elements 
of the multiplicator, whose images in $\otimes^{n-1} U(L/I)$ are 
linearly independent. We shall deal with
several cases. In each of them we shall construct 
elements of $H_2(L/I^n,k)$ that have one parameter $l$, where 
$l\in U(L/I)$. In other words we shall construct a $k$-linear map 
$f:U(L/I)\to H_2(L/I^n,k)\to \otimes^{n-1}U(L/I)$. It is obviously enough
to show that ${\rm ker}f=k\cdot 1$ (since $U(L/I)$ is not finite dimensional). 
In other cases we shall show that ${\rm Im}f$ is not finite dimensional
by proving that it has elements of unbounded degree.

 Recall the elements $\alpha\in I/I'$ and $x\in X$ such that
$a=\phi_x(\alpha)$ was non-zero, and consider all elements of the
form $[\alpha \cdot l,\alpha, \ldots ,\alpha]\otimes 1$, where $l$ is 
any element of $\delta U(L/I)$. Obviously this element is in $I^n$. Its image,
using the mapping $\phi_{x,x,\ldots,x}$ will be $[al,a,\ldots,a]\otimes 1$.
In other words $f(l)=[al,a,\ldots,a]\otimes 1$. Note that if $l\in k\cdot 1$ 
then $f(l)=0$ since in that case $[a\cdot l,a]=0$.
An easy induction shows that 
$$[a,b,b,\ldots,b]\otimes 1 = \sum (-1)^i {{n-1}\choose i}\otimes^ib\otimes a\otimes^{n-1-i}b\otimes 1$$
where $\otimes^i b$ means $b\otimes b\otimes \cdots \otimes b$ ($i$ times). The
referee points out that this formula is known as the Cartan-Weyl formula.
Therefore under the Hopf module isomorphism
$$f(l)=\sum (-1)^i { {n-1}\choose i}(\otimes^ia\otimes al\otimes^{n-2-i}a)\Delta_{n-1}(S(a))$$
$$ + (-1)^{n-1} (\otimes^{n-1}a)\Delta_{n-1}(S(al)).$$
But $S(al) = S(l)S(a)$ so $\Delta_{n-1}(S(al)) = \Delta_{n-1}(S(l)) \Delta_{n-1}(S(a))$ and hence
$$f(l)=[\sum (-1)^i { {n-1}\choose i}(\otimes^ia\otimes al\otimes^{n-2-i}a)$$
$$ + (-1)^{n-1} (\otimes^{n-1}a)\Delta_{n-1}(S(l))]\Delta_{n-1}(S(a)).$$
This can be rewritten as
$$f(l)=(a\otimes a\otimes \cdots \otimes a)[\sum (-1)^i { {n-1}\choose i} (\otimes^i 1 \otimes l \otimes^{n-2-i} 1) +$$
$$ (-1)^{n-1}\Delta_{n-1}(S(l))]\Delta_{n-1}(S(a)).$$

Since $U(L/I)$ is without zero divisors and we are only interested in 
${\rm ker}f$ or the dimension of ${\rm Im}f$, we can consider instead the 
function
$$f(l)=\sum (-1)^i { {n-1}\choose i} (\otimes^i 1 \otimes l \otimes^{n-2-i} 1) + (-1)^{n-1}\Delta_{n-1}(S(l)).$$

In order to compute ${\rm ker}f$, we can apply $\epsilon$ to all but
the $j$-th coordinate of each monomial. This operator, applied to
$\otimes^i 1 \otimes l \otimes^{n-2-i} 1$, yields
$l \delta_{ij}$ ( since $\epsilon(l)=0$ ), while applied to $\Delta_{n-1}(S(l))$ yields
(because $\epsilon$ is a counit) $S(l)$. Therefore for each $0\leq j<n$ 
the result is
$$ (-1)^j { {n-1}\choose j} l + (-1)^{n-1} S(l) = 0.$$
Therefore $S(l) = (-1)^{n+j}{ {n-1}\choose j} l$.

  If $n>2$ we get 
$ S(l) = (-1)^n l$ and $S(l)= (-1)^{n+1} (n-1)l.$

Therefore $(-1)^n l = (-1)^{n+1} (n-1)l$ i.e.
$$ nl = 0 .$$

  As was mentioned above, there are several cases.

  {\bf Case I} If ${\rm char}(k)$ does not divide $n$ and $n>2$ then 
for any $l\in \delta U(L/I)$ we have $f(l)\not =0$ i.e. ${\rm ker}f=k\cdot 1$.

  {\bf Case II} If ${\rm char}(k) \not = 2$. We wish to show
that ${\rm Im}f$ is not finite dimensional. Denoting by $f_1(l)$
the application of
$\epsilon$ to all but the first coordinate, we get
$f_1(l)=l+(-1)^{n-1}S(l)$. This is true also when $n=2$. 
Since $f_1$ is simply $f$ composed with another function, obviously 
${\rm dim(Im}f_1)\leq {\rm dim(Im}f)$. Therefore it is enough to 
consider $f_1$. However, if $x$ is any non-zero Lie 
element in $U(L/I)$ then $S(x^i) = (-1)^i x^i$.
So $f_1(x^i)=x^i+(-1)^{i+n-1}x^i$. Since ${\rm char}(k)\not =2$ then 
for all $i$ of the correct parity we will have $f_1(x^i)=2x^i\not =0$, but
${\rm deg}x^i=i$ will be unbounded, so we are finished.

  {\bf Case III} The only case left is ${\rm char}(k)=2$ and $n$ even.
In this case we still have $f_1(l)=l-S(l)$. Suppose $L/I$ is not 
commutative, therefore there exist $x,y\in L$ such 
that $[x,y]\not \in I$, i.e. $[x,y]\not =0$ in $U(L/I)$.
Consider $l_i=xy^i$. Obviously $S(l_i)=y^ix$, so $f_1(l_i)=xy^i-y^ix=[x,y^i]$.
However the mapping $u\mapsto [x,u]$ is a derivation of $U(L/I)$, and therefore
$$[x,y^i]=\sum_{j=0}^{i-1} y^j[x,y]y^{i-j-1}.$$
Note that $[x,y]y=y[x,y]+[[x,y],y]$, and hence $y^j[x,y]y^{i-j-1}\equiv y^{i-1}[x,y]$
mod $U_{i-1}(L/I)$. Thus $[x,y^i]\equiv iy^{i-1}[x,y]$ mod $U_{i-1}(L/I)$, and
if $i$ is odd then ${\rm deg}f_1(l_i)=i$. Thus the degree of the elements of
the image is unbounded, so the image is infinite dimensional.

  {\bf Case IV} There remains the case where $L/I$ is 
commutative. Thus if $L$ has basis $X$, then $L'\subset I$ so $I/L'\subset L/L'$ is a subspace, and we can
perform a linear change of basis of $L$, so that $I=<L',X_1>$, 
where $X_1$ is a proper subset of $X$. Consider the Lie algebra over 
$\mathbb Z$, $L_1=<Y>$, $I_1=< L_1',Y_1>$, where $Y$ and $Y_1$ are disjoint 
copies of $X$ and $X_1$.
We now use the universal coefficient theorem (see e.g. \cite{Hi} p.176) 
which in our case states that if $k$ is any $\mathbb Z$ module then
$$ 0\to H_2(L_1/I_1^n,\mathbb Z)\otimes_{\mathbb Z}k\to H_2(L_1/I_1^n \otimes_{\mathbb Z}k,k)$$
$$\to Tor_1^{\mathbb Z}(H_1(L_1/I_1^n,\mathbb Z),k)\to 0$$
is exact. Since $H_1(L_1/I_1^n,\mathbb Z)= (L_1/I_1^n)_{ab}=L_1/L_1'$ is a 
free $\mathbb Z$ module then $Tor_1^{\mathbb Z}(H_1(L_1/I_1^n,\mathbb Z),k)=0$.

 Take $k=\mathbb Q$. We have $H_2(L_1/I_1^n,\mathbb Z)\otimes_{\mathbb Z}
\mathbb Q \approx H_2(L_1/I_1^n \otimes_{\mathbb Z}\mathbb Q,\mathbb Q)$. 
However $L_1/I_1^n \otimes_{\mathbb Z}\mathbb Q$ is simply $L_2/I_2^n$ 
where $L_2=<Y>$ and $I_2=<L_2',Y_1>$ taken over $\mathbb Q$. Since 
$\mathbb Q$ has characteristic $0$, we know that $H_2(L_2/I_2^n 
\otimes_{\mathbb Z}\mathbb Q,\mathbb Q)$
is infinite dimensional. Therefore $H_2(L_1/I_1^n,\mathbb Z)$ must also
have infinite torsion-free rank as a $\mathbb Z$-module. Apply now the 
universal coefficient theorem with $k$ any field of characteristic $2$. Again 
$H_2(L_1/I_1^n,\mathbb Z)\otimes_{\mathbb Z}k\approx H_2(L_1/I_1^n 
\otimes_{\mathbb Z}k,k)$. Once again $L_1/I_1^n \otimes_{\mathbb Z}k$ is 
exactly $L/I^n$ of the original Lie algebra. However, since
$H_2(L_1/I_1^n,\mathbb Z)$ has infinite rank then $H_2(L_1/I_1^n,\mathbb Z)
\otimes_{\mathbb Z}k$ is not finitely generated, thus we have proved 
theorem \ref{th:main}.

 Note that in the case $L=<x,y>$, $I=L'$ and $k$ is of characteristic $2$,
even though $H_2(L/I',k)$ is not finitely generated, the image in $U(L/I)$, 
under any of the projections, will be $0$.

\end{document}